\documentclass[11pt,a4paper,leqno]{article}
\usepackage{amsfonts,amssymb,mathrsfs}

\usepackage{amscd}

\vspace{0.5cm}

 4

\font\erm=cmr8

\usepackage{graphicx}

\author{Andrzej Krzysztof Kwa\'sniewski}
\title{A note  on \textit {V}-binomials' recurrence for  V-Lucas sequence companion to U-Lucas sequence}

\newtheorem{defn}{Definition}

\begin{document}

\begin{center}
\noindent { \textbf{A note  on \textit {V}-binomials' recurrence for  V-Lucas sequence companion to U-Lucas sequence}}  \\ 

\vspace{0.5cm} 
Andrzej Krzysztof Kwa\'sniewski

\vspace{0.5cm}
{\erm
	Member of the Institute of Combinatorics and its Applications, 
Winnipeg, Manitoba, Canada \\
	PL-15-674 Bia\l ystok, Konwaliowa 11/11, Poland\\
	e-mail: kwandr@gmail.com\\}
\end{center}

\vspace{1cm}

\noindent \textbf{Summary} 

\noindent Following \cite{Savage} (2009) we deliver \textit{V}-binomials' recurrence formula for Lucas sequence $V = \left\langle V_n\right\rangle_{n\geq 0}$  companion to $ U = \left\langle U_n\right\rangle_{n\geq 0}$ sequence \cite{EdL} (1878). This formula is not present neither in \cite{EdL} (1878)  nor in \cite{Savage} (2009), nor in \cite{Fon}  (1915), nor in \cite{Ward} (1936), nor in \cite{JM} (1949) and neither in all other  quoted here as  "`Lucas $(p,q)$-people"' references [1-34]. Meanwhile  \textit{V}-binomials' recurrence formula for Lucas sequence $V_n$ easily  follows from the original Theorem 17 in \cite{Savage}. 
\noindent Our formula may and should be confronted with \cite{Fon} (1915) Fonten\'e recurrence i.e. (6) or (7) identities in \cite{Gould} (1969) which, as we indicate, also stem easily  from the Theorem 17 in \cite{Savage}.

\vspace{0.2cm}
\noindent AMS Classification Numbers: 05A10 , 05A30

\vspace{0.2cm}

\section{Preliminaries - notation.}

\noindent Notation   a,b $a\neq b$ in \cite{EdL} (1878) is used for the roots of the equation $x^2= P x -Q$  or $(a,b)\equiv(u,v)$ in \cite{Savage} (2009) for the roots of the equation $x^2  =  \ell x -1$.

\noindent The identification  $(a,b) \equiv (p,q)$  i.e.  $p,q$  are  used in "`Lucas $(p,q)$-people"' publications recently and in  recent past (look into not complete list of references [2 -  33]  and  $p,q$-references therein).

\noindent Lucas $(p,q)$-people would then use $U$-identifications:

$$n_{p,q} = \sum_{j=0}^{n-1}{p^{n-j-1}q^j}  = U_n  = \frac{p^n - q^n} {p-q},\  0_{p,q}=U_0 = 0, \ 1_{p,q}=U_1 = 1,$$
where  $p,q$  denote now the roots of the equation $x^2= s x +t \equiv x^2= Px -Q$ hence  $p+ q = s =P$ and $pq = Q = - t$ and the empty sum convention was used for $0_{p,q} = 0$.
\noindent Usually one considers -as we do now- the  $p\neq q$ case. In general also $s\neq t$ - though  according to the context \cite{G-V} (1989) $s=t$  may happen  to be the case of interest.

\vspace{0.2cm}

\noindent The Lucas $U$-binomial coefficients ${n \choose k}_U \equiv {n \choose k}_{p,q}$  are then defined (\cite{EdL} (1878), \cite{Ward} (1936), \cite{JM}(1949), \cite{TF} (1964), \cite{Gould} (1969) ) as follows.

\begin{defn}
	Let  $U$ be as in \cite{EdL} i.e $U_n= n_{p,q}$ then $U$-binomial coefficients for any $n,k \in \mathbb{N}\cup\{0\}$ are defined as follows
	\begin{equation}
		{n \choose k}_U ={n \choose k}_{p,q} = \frac{n_{p,q}!}{k_{p,q}! \cdot (n-k)_{p,q}!} = \frac{n_{p,q}^{\underline{k}}}{k_{p,q}!}
	\end{equation}
	\noindent where $n_{p,q}! = n_{p,q}\cdot(n-1)_{p,q}\cdot ... \cdot 1_{p,q}$ and $n_{p,q}^{\underline{k}} = n_{p,q}\cdot(n-1)_{p,q}\cdot ...\cdot (n-k+1)_{p,q}.$
\end{defn} 

\vspace{0.2cm}

\begin{defn}
Let $V$ be as in \cite{EdL} i.e $V_n = p^n +q^n$, hence $V_0 = 2$ and  $V_n = p +q =s \equiv P$. Then $V$-binomial coefficients for any $n,k \in \mathbb{N}\cup\{0\}$ are defined as follows
	\begin{equation}
		{n \choose k}_V =\frac{V_n!}{V_k!\cdot V_(n-k)!} = \frac{V_n^{\underline{k}}}{V_k!}
	\end{equation}
	\noindent where $V_n! = V_n \cdot V_{n-1}\cdot...\cdot V_1$ and $V_n^{\underline{k}}=V_n \cdot V_{n-1}\cdot ...\cdot V_{n-k+1}.$
\end{defn}

\vspace{0.2cm}

\noindent One easily generalizes $L$-binomial to  $L$-multinomial coefficients \cite{MD} .

\begin{defn}\label{def:symbol}
Let $L$ be any natural numbers' valued sequence i.e. $L_n\in\mathbb{N}$ and $s\in\mathbb{N}$. \textbf{$L$-multinomial coefficient} is then identified with the symbol

\begin{equation}
	{n \choose {k_1,k_2,...,k_s}}_L = \frac{L_n!}{L_{k_1}!\cdot ... \cdot L_{k_s}!}
\end{equation}

\vspace{0.2cm}
\noindent where $k_i\in\mathbb{N}$ and $\sum_{i=1}^{s}{k_i} = n$ for $i=1,2,...,s$. Otherwise it is equal to zero.
\end{defn}

\noindent Naturally  for any natural $n,k$ and $k_1+...+k_m=n-k$ the following holds

\begin{equation} \label{eq:mult1}
{n \choose k}_L \cdot {n-k \choose {k_1,k_2,...,k_m}}_L  = {n \choose {k,k_1,k_2,...,k_m}}_L 
\end{equation}


\section{$V$-binomial coefficients' recurrence}

The authors of \cite{Savage}   prove  trivially an observation named the Theorem 17 i.e  the following nontrivial recurrence for the general case of ${r+s \choose r,s}_{L[p,q]}$ $L$-binomial arrays in multinomial notation.
Let $s,r>0$. Let $L = L[p,q]$ be any zero characteristic field valued sequence ($L_n \neq 0$). Then

\begin{equation}
{r+s \choose r,s}_{L[p,q]}=g_1(r,s)\cdot{r+s-1 \choose r-1,s}_{L[p,q]}+g_2(r,s)\cdot{r+s-1 \choose r,s-1}_{L[p,q]}
\end{equation}
where   $ {r \choose r,0}_L = {s \choose 0,s}_L =1.$

\begin{equation}
 L[p,q]_{r+s} =  g_1(r,s) \cdot L[p,q]_r   +  g_2(r,s) \cdot L[p,q]_s.                               
\end{equation}

\vspace{0.3cm}

\noindent \textbf{Compare} the above now  common knowledge formulas (5) and (6) \textbf{with }: [JM] formula (2) and formula in between (1') and (2) in \cite{JM} (1949) or [Gould] compare formula (7) in \cite{Gould} (1969) with  this note formula(5), or [KaKi] see formulas (51) and (40) in \cite{KaKi} (1992) or [G-V] see this correspondence in section (10.1) of  \cite{G-V} (1989)  
or [K-W] see this correspondence in \cite{K-W} (1989)  or [Corcino] compare this note formula (5) with the Theorem 1  formulas (13) and (14) - with the simple proof just checking - in \cite{Corsino} (2008)  or [MD1] compare this note formulas (5) and (6) with  \textit{the special case} formulas (2) and (1) in \cite{MD2} \textbf{v[1]} [here only  special form $T_\lambda$ tiling cobweb admissible natural numbers' valued sequences are admitted] where note that (1) and (2) formulas are given original combinatorial interpretations in terms of tilings of the so called cobweb posets  and (2) is given a combinatorial proof  or [MD2] compare this note formulas (5) and (6) with (16) and (11) in  \cite{MD2} \textbf{v[2]} where a trivial derivation of nontrivial (16) is supplied - in reverse order with respect to the corresponding derivation in \cite{Savage}.

\vspace{0.3cm}

\noindent \textbf{Historical Note}
\noindent As accurately noticed by Knuth and Wilf in \cite{K-W}  the recurrent relations for Fibonomial coefficients appeared already in 1878  Lukas work (see: p. 27 , formula (58) in \cite{EdL}. In our opinion - Lucas  Th\'eorie des Fonctions Num\'eriques Simplement P\'eriodiques is the far more non-accidental context  for binomial-type coefficients exhibiting their relevance to hyperbolic trigonometry  [see for example \cite{ Bakk}, \cite{Bajguz}].  

\noindent Indeed. Taking here into account the  $U$-addition formula i.e. the first of two trigonometric-like $L$-addition formulas (42) from \cite{EdL} [see also \cite{ Bakk}, \cite{Bajguz}]          ($L[p,q]  = L = U,V$)  i.e.

\begin{equation}
2 U_{r+s} =  U_r V_s  +  U_s V_r ,\ \ \ \ 
2 V_{r+s} =  V_r V_s  +  U_s U_r 
\end{equation}
one readily recognizes that  the $U$-binomial recurrence from the Corollary 18 in  \cite{Savage} is identical with the $U$-binomial recurrence (58) \cite{EdL}.
\noindent See also Proposition 2.2. in \cite{ BSCS}  (2010).

\vspace{0.2cm}

\noindent However  there is no companion  $V$-binomial recurrence neither in  \cite{EdL} (1878) nor in \cite{Savage} (2009)as well as all other quoted papers. Here let us note
on the way as aside remark, that we know \cite{Eric} quite promising analogues of addition rules for both companion sequences. Namely 

\vspace{0.2cm}

\begin{equation}
U_{r+s} =  U_r V_s  - p^n q^n U_{r-s},\ \ \ \ 
V_{r+s} =  V_r V_s  - p^n q^n V_{r-s}. 
\end{equation}

\vspace{0.3cm}

\noindent Now, it is not difficult to realize that the looked for  $V$-binomial recurrence may be given in the form of (5) adapted to  $ L[p,q] = V[p,q]= V $ - Lucas sequence case. Here it is.

\begin{equation}
{r+s \choose r,s}_{V[p,q]} = h_1(r,s){r+s-1 \choose r-1,s}_{V[p,q]} + h_2(r,s){r+s-1 \choose r,s-1}_{V[p,q]},
\end{equation}
where $p \neq q$ and  ${r \choose r,0}_L = {s \choose 0,s}_L =1,$  and  for  $ r \neq s $

\begin{equation}
V_{r+s} =  h_1(r,s)V_r   +  h_2(r,s) V_s,  \  
V_{2r} =  (h_1(r,r) +  h_2(r,r))  \cdot  V_r
\end{equation}
and where  ($p \neq q$) and  $r \neq s $  while 

\begin{equation}
h_1 \cdot (p^r q^s  -  q^r p^s ) = p^{r+s} q^s  -  q^{r+s} p^s,  \ r \neq s , \
\end{equation}
and

\begin{equation}
h_2 \cdot (q^r p^s  -  p^r q^s ) = p^{r+s} q^r  -  q^{r+s} p^r   \ r \neq s , \ 
\end{equation}
Now for $ r = s$ -  having in mind that  $V_k = p^k + q^k$  - a taking into account $p$ interchange with $q$ symmetry - we may establish the identification  (13) for $P \neq 0$.
\noindent (Recall \cite{EdL} (1878), that $p,q$ are roots of $x^2 - Px +Q$ therefore  $p+q = P$  and  $p\cdot q = Q$.)

\begin{equation}
h_1(r,r)= \frac{p^{2r}} {p^r + q^r} ,  \    h_2(r,r) = \frac{q^{2r}} {p^r + q^r} .
\end{equation}
The matters are much easier in $L_n = U[p,q]_n = n_{p,q} = \frac{p^n - q^n} {p-q}$ $U$-Lucas sequence well elaborated  case. See for example formula (2) and formula in between (1') and (2) in \cite{JM} (1949) and come back to acapit starting with  \textbf{Compare} - above. 

\vspace{0.2cm}

\noindent Indeed. One may proceed as above with  $V$-Lucas sequence immediately noting that for $r=s$   $h_1 + h_2 = p^r + q^r$ as it should be, See below.
\vspace{0.1cm} 

\noindent  The recurrent relations (13) and (14) in \cite{Corsino}  for $n_{p,q}$-binomial coefficients are special cases of this paper formula (5) i.e. of Th. 17 in \cite{Savage}
with straightforward identifications of $g_1, g_2$  in (13)  and  in (14) in \cite{Corsino} as well as this paper recurrence (6) for $L = U[p,q]_n = n_{p,q}$ sequence.

\begin{equation}
 g_1 = p^r ,\ \  g_2 =  q^s, 
\end{equation}
or
\begin{equation}
 g_1 = q^r ,\ \  g_2 =  p^s,
\end{equation}
while

\begin{equation}
(s + r)_{p,q} = p^s r_{p,q} +  q^r s_{p,q} = (r + s)_{q,p} = q^r s_{p,q} +  p^s  r_{p,q}.
\end{equation}
hence - in multinomial notation and choosing (14) we have

\begin{equation}
{r+s \choose r,s}_U \equiv {r+s \choose r,s}_{p,q}= p^r\cdot{r+s-1 \choose r-1,s}_U + q^s \cdot{r+s-1 \choose r,s-1}_U,
\end{equation}
where   $ {r \choose r,0}_U = {s \choose 0,s}_U =1.$

\vspace{0.4cm}
\noindent  Now let $A$ be any natural numbers' or even complex numbers' valued sequence. One readily sees that also (1915) Fonten\'e recurrence for Fonten\'e-Ward generalized $A$-binomial coefficients i.e. equivalent identities (6) , (7) in \cite{Gould} \textbf{are special cases of} this paper formula \textbf{(5)} i.e. of Th. 17 in \cite{Savage} with straightforward identifications of $g_1, g_2$  in this paper formula (5)  while this paper recurrence  (6) becomes trivial identity.

\noindent Namely, the identities (6) and (7) from \cite{Gould} (1969) read correspondingly:

\begin{equation}
{r+s \choose r,s}_A = 1 \cdot {r+s-1 \choose r-1,s}_A + \frac{A_{r+s} - A_r}{A_s}{r+s-1 \choose r,s-1}_A,
\end{equation}

\begin{equation}
{r+s \choose r,s}_A = \frac{A_{r+s} - A_s} {A_{r}} \cdot {r+s-1 \choose r-1,s}_A +  1 \cdot {r+s-1 \choose r,s-1}_A,
\end{equation}
where $p \neq q$ and  $ {r \choose r,0}_L = {s \choose 0,s}_L =1.$ And finally we have tautology identity

\begin{equation}
A_{s + r} \equiv \frac{A_{r+s} - A_s}{A_r}\cdot A_r +  1 \cdot A_s.
\end{equation}

\vspace{0.3cm}

\noindent As for \textbf{combinatorial interpretations} of $L$-binomial or $L$-multinomial coefficients we leave that subject apart from this note 
because this note is to be deliberately short. Nevertheless we direct the reader to some papers and references therein; these are here the following: \cite{ BSCS} (2010),\cite{Savage} (2009), 
\cite{MD} (2009), \cite{MD2} (2009 , \cite{G-V} (1989), \cite{wachs} (1991),\cite{KaKi} (1992),  \cite{medicis} (1993),  \cite{wachs 2} (1994), 
\cite{Park} (1994),  \cite{RW} (2004), and to this end see \cite{Voigt}. The list is far not complete. For example the on combinatorial interpretation work of Arthur T. Benjamin $http://www.math.hmc.edu/~benjamin/$ should  be notified also.

\vspace{0.4cm}

\noindent \textbf{Final remark}. The above presentation, definitions and recurrent formulas like  (5)  and (6) extend  correspondingly to Horadam $W$  and $H$ sequences while (10)-(12) just stay the same  under the replacement: $V \rightarrow H $  i.e     $V_n = p^n +q^n  \rightarrow H_n = A \cdot p^n +  B \cdot q^n$  while  (13) becomes

\begin{equation}
h_1(r,r)= \frac{A \cdot p^{2r}} {A \cdot p^r + B \cdot q^r} ,  \    h_2(r,r) = \frac{B \cdot q^{2r}} {A \cdot p^r + B \cdot q^r}.
\end{equation}

\vspace{0.2cm}

\noindent To this end recall that Lucas sequences are special cases of Horadam sequences i.e. the extension looks like: (Lucas \cite{EdL} (1878)) $U_n =  \frac{p^n - q^n} {p-q} \rightarrow   \frac{A \cdot p^n - B\cdot q^n} {p-q} = W_n  $ (Horadam : \cite{Horadam1}, \cite{Horadam2}  (1965), \cite{Horadam3}, \cite{Elmore}, \cite{Horadam4}.  

\vspace{0.2cm}

\noindent To this end note also that Horadam $W$ sequence is a special case of Horadam H sequence  with obvious choice of  $A$  and  $B$ in the latter.

\vspace{0.2cm}

\noindent A more detailed presentation of Horadam binomials' properties we leave to the subsequent note.



\end{document}